\documentclass[11pt,reqno]{amsart}

\usepackage{color, MnSymbol, bm}
\usepackage{multirow}

\usepackage{natbib}
\usepackage{hyperref}
\usepackage[latin1]{inputenc}
\usepackage{graphicx}
\usepackage{bbm}
\usepackage{subcaption}
\usepackage[autostyle]{csquotes}
\usepackage{siunitx}
\usepackage{float}
\numberwithin{equation}{section}

\theoremstyle{plain}                    % default
\newtheorem{thm}{Theorem}[section]
\newtheorem{lem}{Lemma}[section]

\newtheorem{prop}{Proposition}[section]

\theoremstyle{definition}

\theoremstyle{remark}
\newtheorem{rem}{Remark}[section]               %\renewcommand{\therem}{}

\def\N{{\mathbb N}}
\def\Z{{\mathbb Z}}

\def\P{{\mathbb P}}
\def\E{{\mathbb E}}
\def\F{{\mathcal F}}
\def\1{{\mathbbm{1}}}

\newcommand{\var}{\mathop{\rm var}\nolimits}

\newcommand{\Pois}{\mathop{\rm Pois}\nolimits}
\newcommand{\Skellam}{\mathop{\rm Skellam}\nolimits}
\newcommand{\Bin}{\mathop{\rm Bin}\nolimits}

\newcommand{\be}{\begin{equation}}
\newcommand{\bd}{\begin{displaymath}}
\newcommand{\ed}{\end{displaymath}}
\newcommand{\bea}{\begin{eqnarray}}
\newcommand{\eea}{\end{eqnarray}}
\newcommand{\bean}{\begin{eqnarray*}}
\newcommand{\eean}{\end{eqnarray*}}

\begin{document}

\thispagestyle{empty}

\baselineskip13pt

\begin{center}
{\large \sc Mixing properties of integer-valued GARCH processes}
\end{center}

\vspace*{1cm}

\begin{center}
Paul Doukhan\\
Universit\'e Cergy-Pontoise\\
UMR 8088 Analyse, G\'eom\'etrie et Mod\'elisation\\
2, avenue Adolphe Chauvin\\
95302 Cergy-Pontoise Cedex\\
France\\
E-mail: doukhan@cyu.fr\\[1cm]
Naushad Mamode Khan\\
University of Mauritius\\
Department of Economics and Statistics\\
Reduit 80837\\
Mauritius\\
E-mail: n.mamodekhan@uom.ac.mu\\[1cm]
Michael H.~Neumann\\
Friedrich-Schiller-Universit\"at Jena\\
Institut f\"ur Mathematik\\
Ernst-Abbe-Platz 2\\
D -- 07743 Jena\\
Germany\\
E-mail: michael.neumann@uni-jena.de\\[1.5cm]
\end{center}

%\vspace*{1.0cm}

\begin{center}
{\bf Abstract}
\end{center}
We consider models for count variables with a GARCH-type structure.
Such a process consists of an integer-valued component and a volatility process.
Using arguments for contractive Markov chains we prove that this bivariate process has a unique
stationary regime. Furthermore, we show absolute regularity ($\beta$-mixing) with
geometrically decaying coefficients for the count process.
These probabilistic results are complemented by a statistical analysis and
a few simulations.
\vspace*{2cm}

\footnoterule \noindent {\sl 2010 Mathematics Subject
Classification:} Primary 60G10; secondary 60J05. 
{\sl Short title:} Mixing of INGARCH processes. \vfill
\noindent
version: \today

\newpage

\setcounter{page}{1}
\pagestyle{headings}
\normalsize

%%%%%%%%%%%%%%%%%%%%%%%%%%%%%%%%%%%%%%%%%%%%%%%%%%%%%%%%%%%%%%%%%%%%%%%%%%%%%%%
\section{Introduction and notation}
\label{S1}
%%%%%%%%%%%%%%%%%%%%%%%%%%%%%%%%%%%%%%%%%%%%%%%%%%%%%%%%%%%%%%%%%%%%%%%%%%%%%%%

Models involving integer-valued random variables have attracted increasing attention
in the recent years. In most cases, the random variables are assumed to be
non-negative, since they represent numbers of counts as e.g. the series of road traffic accidents,
the monthly or annual numbers of deaths due to accidents or diseases, or more recently,
the striking daily numbers of infected people and deaths due the novel corona virus.
Sometimes, it is however necessary to allow for both non-negative and negative integer-valued
random variables. A typical field of application is the description of score
differences in sports, e.g.~the number of goals of the home team minus that
of the away team. Another common example is that of price changes in finance,
such as the tick by tick data in a Trade and Quote database, that are often
represented as a combination of positive and negative integer values.
Furthermore, similar types of observations also appear in differenced series
that are initially non-stationary. While in many applications a serial dependence between the observations is not
taken into account, we consider models which allow to describe and exploit
dependencies between consecutive variables.
In view of their popularity in financial time series analysis and because of
their flexibility we focus here on models with a GARCH-type structure.
Adapting the structure of the classical GARCH model by \citet{Bol86},
\citet{FT11,FT12} considered such models for non-negative count variables
and investigated conditional maximum likelihood estimators of the corresponding parameters.
These authors assumed that the count variable~$X_t$ at time~$t$ conditioned on the past
has a Poisson distribution with an intensity~$\lambda_t$ which itself is random
and depends on lagged values of the count and intensity processes.
Since in case of a Poisson distribution the variance is equal to the mean, a GARCH-type
structure is imposed by the equation
\begin{equation}
\label{1.1}
\lambda_t \,=\, \omega \,+\, \alpha_1 X_{t-1} \,+\, \cdots \,+\, \alpha_p X_{t-p}
\,+\, \beta_1 \lambda_{t-1} \,+\, \cdots \,+\, \beta_q \lambda_{t-q}
\end{equation}
or by nonlinear variants, $\lambda_t=f_\theta\big(X_{t-1},\ldots,X_{t-q},\lambda_{t-1},\ldots,\lambda_{t-q}\big)$,
where $\theta$ is a suitable parameter.

In this paper, we consider once more processes with a GARCH-type structure. In contrast to
the papers mentioned above, we allow for integer-valued variables which can attain both
non-negative and negative values. The most prominent example is the distribution
introduced by \citet{Ske46}, which is the distribution of the difference of two
independent Poisson variates with respective parameters $\lambda_1$ and $\lambda_2$,
In the special case of $\lambda_1=\lambda_2$ considered by \citet{Irw37}, the corresponding distribution has zero mean.
Therefore and in contrast to (\ref{1.1}), the conditional mean is no longer suitable
to generate a GARCH-type structure.
We will focus on second moments and consider processes where the integer-valued variables~$X_t$
have a conditional distribution $Q_{v_t}$ where $\int x^2\,dQ_{v}(x)=v$ and
\begin{equation}
\label{1.2}
v_t \,=\, f\big(X_{t-1}^2,\ldots,X_{t-p}^2,v_{t-1},\ldots,v_{t-q}\big).
\end{equation}
\citet{AAO18} considered such a Skellam-GARCH process of order $p=q=1$ and derived the
estimating equations for a conditional maximum likelihood estimator of the parameters.
However, perhaps because of the absence of suitable probabilistic tools for such models,
they did not provide a further analysis of the asymptotic properties of this estimator.
These authors also provided an overview of related results and applied the model
to differences of non-negative data of counts of monthly drug crimes.

In this contribution, we primarily focus on stochastic properties such as existence and uniqueness
of a stationary distribution and absolute regularity of integer-valued GARCH processes.
In the related case of Poisson-GARCH processes with linear or nonlinear specifications
$\lambda_t=f\big(X_{t-1},\ldots,X_{t-p},\lambda_{t-1},\ldots,\lambda_{t-q}\big)$, there are
several forerunners of the present work and it turns out that we can build on the methods
derived there.
Mixing properties of such processes have been derived for a first time in \citet{Neu11},
for INGARCH(1,1) processes under a contractive condition.
This has been generalized in \citet{DN19} for the INGARCH($p$,$q$) case and
under a weaker semi-contractive condition which resulted in a somewhat unusual subexponential
decay of the mixing coefficients. \citet{DLN20} proved absolute regularity of the count process
again in the INGARCH(1,1) case but allowing a possibly non-stationary (explosive)
behavior of the process. Finally, \citet{Neu20} proved absolute regularity with an
exponential decay of the mixing coefficients for INGARCH($p$,$q$) processes
under a fully contractive condition,
\begin{equation}
\label{1.3}
\big| f(x_1,\ldots,x_p,\lambda_1,\ldots,\lambda_q) \,-\, f(x_1',\ldots,x_p',\lambda_1',\ldots,\lambda_q') \big|
\,\leq\, \sum_{i=1}^p c_i |x_i-x_i'| \,+\, \sum_{j=1}^q d_j |\lambda_j-\lambda_j'|,
\end{equation}
where $c_1,\ldots,c_p,d_1,\ldots,d_q$ are non-negative constants are such that $\sum_{i=1}^p c_i+\sum_{j=1}^q d_j<1$.
We will also impose the contractive condition (\ref{1.3}) on the volatility function~$f$,
however, in contrast to the papers mentioned above, the arguments of this function
reflect second-order properties of the process $(X_t)_t$.
Note that $(Y_t)_t$ and $(Z_t)_t$ with $Y_t=\big(X_t,\ldots,X_{t-p+1},v_t,\ldots,v_{t-q+1}\big)$
and $Z_t=\big(X_t^2,\ldots,X_{t-p+1}^2,v_t,\ldots,v_{t-q+1}\big)$
are both (first-oder) Markov chains.
We show in Section~\ref{S2} that the contractive condition on~$f$ yields a contraction property for $(Z_t)_t$
in terms of a suitable Wasserstein metric.
This implies by the Banach fixed point theorem that $(Z_t)_t$ possesses
a unique stationary distribution, and a simple extra argument shows that the same property
holds true for the process $(Y_t)_t$ which is of actual interest here.
Furthermore, we use the contraction property once more to prove almost effortlessly absolute regularity
($\beta$-mixing) with exponentially decaying coefficients of the count process $(X_t)_t$.

We are convinced that these results can serve as a basis for further work with such models
without any hassle.
As an example, population dynamics can be considered after differentiation,
in order to rate the speed or the acceleration of the evolution of species under consideration;
indeed both characteristics may be either positive or negative.
As an illustration of their usefulness, we apply in Section~\ref{S3} our results to prove
asymptotic normality of a least squares estimator of the parameters of a Skellam-ARCH model.
All proofs and a few auxiliary results are collected in a final Section~\ref{S4}.
\bigskip

%%%%%%%%%%%%%%%%%%%%%%%%%%%%%%%%%%%%%%%%%%%%%%%%%%%%%%%%%%%%%%%%%%%%%%%%%%%%%%%
\section{Main results}
\label{S2}
\subsection{Assumptions and a preview of the results}
\label{SS2.1}
%%%%%%%%%%%%%%%%%%%%%%%%%%%%%%%%%%%%%%%%%%%%%%%%%%%%%%%%%%%%%%%%%%%%%%%%%%%%%%%

We consider a class of integer-valued processes $(X_t)_{t\in\Z}$ defined on some probability
space $(\Omega,\F,\P)$, where, for all $t\in\Z$,
\begin{subequations}
\begin{eqnarray}
\label{21.1a}
X_t\mid \F_{t-1} & \sim & Q_{v_t}, \\
\label{21.1b}
v_t & = & f(X_{t-1}^2,\ldots,X_{t-p}^2,v_{t-1},\ldots,v_{t-q}),
\end{eqnarray}
and $\F_s=\sigma(X_s,\lambda_s,X_{s-1},\lambda_{s-1},\ldots)$ denotes the $\sigma$-field
generated by the random variables up to time~$s$.
Assuming that~$f$ takes values in some set $V\subseteq [0,\infty)$, this function has to be defined
on $\N_0^p\times V^q$.
The parameter~$v_t$ stands for the conditional second moment of~$X_t$, i.e.~the family of distributions $(Q_v)_{v\in V}$
on $(\Z, {\mathcal P}(\Z))$ is parametrized such that
\begin{equation}
\label{21.1c}
\int x^2 \, dQ_v(x) \,=\, v \qquad \forall v\in V.
\end{equation}
\end{subequations}
A frequently considered special case is that of a linear model of order~$p$, $q$, where
\begin{equation}\label{LGARCH}
f(x_1,\ldots,x_p,v_1,\ldots,v_q)
\,=\, \omega \,+\, \sum_{i=1}^p \alpha_i x_i \,+\, \sum_{j=1}^q \beta_j v_j.
\end{equation}
It is clear that the processes $\bm{Y}=(Y_t)_{t\in\Z}$
and $\bm{Z}=(Z_t)_{t\in\Z}$ with
$Y_t=(X_t,\ldots,X_{t-p+1},v_t,\ldots,v_{t-q+1})$
and $Z_t=(X_t^2,\ldots,X_{t-p+1}^2,v_t,\ldots,v_{t-q+1})$ are
time-homogeneous Markov chains with state spaces $S=\Z^p\times V^q$ and $S^\geq =\N_0^p\times V^q$,
respectively.
The following conditions ensure existence and uniqueness of a stationary distribution
of $(Z_t)_{t\in\Z}$, and eventually of $(Y_t)_{t\in\Z}$ as well.
Furthermore, they also yield absolute regularity of the count process $(X_t)_{t\in\Z}$.
\bigskip
\begin{itemize}
\item[{\bf (A1)}]
There exist non-negative constants $c_1,\ldots,c_p,d_1,\ldots,d_q$ such that
$\displaystyle \sum_{i=1}^p c_i+\sum_{j=1}^q d_j<1$ and
\begin{displaymath}
\left| f(x_1,\ldots,x_p,v_1,\ldots,v_q) \,-\, f(x_1',\ldots,x_p',v_1',\ldots,v_q') \right|
\,\leq\,  \sum_{i=1}^p c_i\; |x_i - x_i'| \,+\, \sum_{j=1}^q d_j\; |v_j - v_j'|
\end{displaymath}
holds for all $(x_1,\ldots,x_p,\lambda_1,\ldots,\lambda_q)$, $(x_1',\ldots,x_p',\lambda_1',\ldots,\lambda_q')\in S^\geq$.\\
\item[{\bf (A2)}]
The family of distributions $(Q_v)_{v\in V}$ is increasing in the following sense:
If, $v<v'$, $X\sim Q_v$, and $X'\sim Q_{v'}$, then~$|X|$ is stochastically not greater than~$|X'|$, i.e.
\begin{displaymath}
P\big(|X|\leq k\big) \,\geq\, P\big( |X'|\leq k \big), \qquad \forall k\in\N.
\end{displaymath}
\end{itemize}
\bigskip

\noindent
{\bf Examples}
\begin{itemize}
\item[1)] {\bf Symmetric Skellam distributions}\\
Let, for $v\in V$, $Q_v=\Skellam(v/2,v/2)$ be a Skellam distribution with parameters $v/2$, $v/2$, i.e.
$Q_v$ is the distribution of two independent Poisson variates with parameter $v/2$ each.
Suppose that $X_1,X_2\sim\Pois(v/2)$ are independent. Then
\begin{displaymath}
\int x^2\, dQ_v(x) \,=\, \var(X_1-X_2) \,=\, \var(X_1)+\var(X_2) \,=\, v,
\end{displaymath}
i.e.~(\ref{21.1c}) is satisfied.
To prove that {\bf (A2)} is fulfilled by the family $(Q_v)_{v\in V}$, suppose that $v<v'$.
If $Y\sim\Skellam(v/2,v/2)$ and $Z\sim\Skellam((v'-v)/2,(v'-v)/2)$ are independent,
it follows from the properties of Poisson distributions  that $Y+Z\sim\Skellam(v'/2,v'/2)$.
Since the probability mass function of~$Y$ is symmetric and unimodal
(see e.g. \citet{AO10}) we have that
\bd
P( |Y| \leq k ) \,\geq\, P( |Y+l| \leq k ), \qquad \forall (l,k)\in\Z\times\N_0,
\ed
which implies that
\bd
P\left( |Y| \leq k \right) \,\geq\, \sum_{l\in\Z} P( |Y+l| \leq k ) P( Z=l )
\,=\, P( |Y+Z| \leq k ), \qquad \forall k\in\N_0.
\ed
Hence, $|Y|$ is stochastically not greater than $|Y+Z|$.
\item[2)] {\bf Mixtures of symmetric Skellam distributions}\\
Let $G$ be the distribution of a non-negative random variable with
$\int_{[0,\infty)} x\, dG(x)=\mu\in(0,\infty)$.

Then
\begin{displaymath}
Q_v \,=\, \int_{[0,\infty)} \Skellam\big( \frac{vx}{2\mu}, \frac{vx}{2\mu} \big) \, dG(x)
\end{displaymath}
is a mixture of symmetric Skellam distributions.
We have that
\begin{displaymath}
\int x^2\, dQ_v(x) \,=\, \int_{[0,\infty)} \int x^2 \,d\Skellam\big( \frac{vx}{2\mu}, \frac{vx}{2\mu} \big) \, dG(x)
\,=\, \int_{[0,\infty)} \frac{vx}{\mu} \, dG(x) \,=\, v,
\end{displaymath}
and, for $v<v'$,
\begin{eqnarray*}
Q_v\big( \{-k,\ldots,k\} \big)
& = & \int_{[0,\infty)} \Skellam\big( \frac{vx}{2\mu}, \frac{vx}{2\mu} \big)(\{-k,\ldots,k\}) \, dG(x) \\
& \geq & \int_{[0,\infty)} \Skellam\big( \frac{v'x}{2\mu}, \frac{v'x}{2\mu} \big)(\{-k,\ldots,k\}) \, dG(x) \\
& = & Q_{v'}\big( \{-k,\ldots,k\} \big),
\end{eqnarray*}
i.e.~(\ref{21.1c}) and {\bf (A2)} are satisfied.
\\
A notable special case is that of a zero-inflated Skellam distribution,
where $B\sim G$ follows a Bernoulli distribution with parameter $p\in(0,1)$.
If $B$ and $X\sim\Skellam(v/(2p),v/(2p))$
are independent, then $BX$ has a zero-inflated Skellam distribution.
Such a distribution was used by \citet{KN06} and \citet{AK14} to account
for an excess of zero counts in certain medical data.
\item[3)] {\bf Poisson distributions}\\
If $X\sim\Pois(\lambda)$, then $EX^2=\lambda^2+\lambda$. This is equal to $v>0$
if and only if $\lambda=\sqrt{v+1/4}-1/2$. In order to obey (\ref{21.1c})
we choose $Q_v=\Pois(\sqrt{v+1/4}-1/2)$.
For $v<v'$ we have that $\sqrt{v+1/4}-1/2<\sqrt{v'+1/4}-1/2$. Since $X\sim\mbox{Pois}(\sqrt{v+1/4}-1/2)$
is stochastically not greater than $X\sim\Pois(\sqrt{v'+1/4}-1/2)$ we see that condition
{\bf (A2)} is satisfied.
\item[4)] {\bf Mixtures of Poisson distributions}\\
Let, as in Example~2, $G$ be the distribution function of a non-negative random variable with
$\int_0^\infty x\, dG(x)=\mu\in(0,\infty)$.
Then
\begin{displaymath}
Q_v \,=\, \int_0^\infty \Pois\big( \sqrt{vx/\mu + 1/4} - 1/2 \big) \, dG(x)
\end{displaymath}
is a mixture of Poisson distributions. Then condition {\bf (A2)} is obviously fulfilled.
Furthermore, since
\begin{displaymath}
\int x^2\, dQ_v(x) \,=\, \int_{[0,\infty)} \int x^2\, d\Pois\big( \sqrt{vx/\mu + 1/4} - 1/2 \big) \, dG(x)
\,=\, \int_{[0,\infty)} \frac{vx}{\mu} \, dG(x) \,=\, v
\end{displaymath}
we see that (\ref{21.1c}) is also satisfied.
Poisson distributions can be used for modeling data from various fields, e.g.~the number of
financial transactions within a certain time period or the number of claims in an insurance context.
When dealing with a collection of individual transactions corresponding to different trading
strategies or with a collection of claim numbers from persons with different features (age, health state,...)
an appropriate mixture of Poisson distributions seems to be more adequate.
Notable special cases are that of a zero-inflated Poisson distribution which appears
in case of $G$ a Bernoulli distribution  with parameter $p\in (0,1)$ or a negative binomial distribution,
if $G$ has a Gamma distribution. %; see e.g.~... \textcolor[rgb]{1,0,0}{\bf Reference?}
A negative binomial distribution is often preferred to a Poisson distribution
if data are overdispersed, i.e.~if their variance is greater than their mean (as this is the case for all mixed Poison distributions).
\item[5)] {\bf Binomial distributions}\\
A $\Bin(n,p)$ distribution ($p\in(0,1)$) can be used for modeling underdispersed data since its variance~$np(1-p)$
is less than its mean~$np$. To satisfy (\ref{21.1c}), we set
\begin{displaymath}
Q_v \,=\, \Bin\big( n, g(v) \big),
\end{displaymath}
where $g\colon\; (0,n^2)\rightarrow (0,1)$ is a strictly monotonic function such that
\begin{displaymath}
\int x^2\, d\Bin(n, g(v)) \,=\, ng(v) \,+\, n(n-1)\big( g(v)\big)^2 \,=\, v.
\end{displaymath}
To see that {\bf (A2)} is fulfilled, let $U_1,\ldots,U_n$ be independent and uniformly distributed on $[0,1]$.
Let $v,v'\in(0,n^2)$. Then $X:=\sum_{i=1}^n\1_{\{U_i\le g(v)\}}\sim Q_v$ and $X':=\sum_{i=1}^n\1_{\{U_i\le g(v')\}}\sim Q_{v'}$.
If $v<v'$, then it follows from the construction that $X\leq X'$ with probability one which implies
that~$X$ is stochastically not greater than~$X'$.
\item[6)] {\bf Some asymmetric distributions over $\Z$}\\
Let $Q_v$ ba any of the above distributions and let $Y\sim Q_v$ and $R$ be independent, where $P(R\in\{-1,1\})=1$.
If $P(R=1)\neq 1$, then the distributions of the random variable $X=RY$ is asymmetric over~$\Z$ and obeys (\ref{21.1c}).
The corresponding family of distributions satisfies {\bf A2}) as soon as $\big(Q_v\big)_{v\in V}$ does.
\end{itemize}
\bigskip

\noindent
In the following we derive a contraction property of $(Z_t)_t$ in terms of
a suitable Wasserstein metric.
As shown in \citet[Chapter~3]{Ebe17} and \citet[Theorem 20.3.4]{DMPS18},
this implies by the Banach fixed point theorem that $(Z_t)_t$ possesses
a unique stationary distribution. A simple extra argument shows that the same property
holds true for the process $(Y_t)_t$ which is of actual interest here.
Furthermore, we use the contraction property once more to prove almost effortlessly absolute regularity
($\beta$-mixing) with exponentially decaying coefficients of the count process $(X_t)_t$.
\bigskip

%%%%%%%%%%%%%%%%%%%%%%%%%%%%%%%%%%%%%%%%%%%%%%%%%%%%%%%%%%%%%%%%%%%%%%%%%%%%%%%
\subsection{Contraction}
\label{SS2.2}
%%%%%%%%%%%%%%%%%%%%%%%%%%%%%%%%%%%%%%%%%%%%%%%%%%%%%%%%%%%%%%%%%%%%%%%%%%%%%%%

First of all, we transfer the contraction condition {\bf (A1)} for the intensity process
into a contraction property for the~$Z_t$.\medskip

We consider the following metric on~$S^\geq$:
\begin{displaymath}
\Delta_{\gamma,\delta}\big( (x_1,\ldots,x_p,v_1,\ldots,v_q), (x_1',\ldots,x_p',v_1',\ldots,v_q') \big)
\,=\, \sum_{i=1}^p \gamma_i\; |x_i-x_i'| \,+\, \sum_{j=1}^q \delta_j\; |v_j-v_j'|,
\end{displaymath}
where $\gamma_1,\ldots,\gamma_p,\delta_1,\ldots,\delta_q$ are strictly positive constants.
Let $y=(x_1,\ldots,x_p,v_1,\ldots,v_q)$, $y'=(x_1',\ldots,x_p',v_1',\ldots,v_q')\in S$
be arbitrary and, accordingly
$z=(x_1^2,\ldots,x_p^2,v_1,\ldots,v_q)$, $z'=({x_1'}^2,\ldots,{x_p'}^2,v_1',\ldots,v_q')\in S^\geq$.
With an appropriate choice of $\gamma_1,\ldots,\gamma_p$, $\delta_1,\ldots,\delta_q$,
we can construct random vectors
$Y=(X,x_1,\ldots,x_{p-1},\lambda,\lambda_1,\ldots,\lambda_{q-1})$ and\\
$Y'=(X',x_1',\ldots,x_{p-1}',\lambda',\lambda_1',\ldots,\lambda_{q-1}')$ on a suitable
probability space $(\widetilde{\Omega},\widetilde{\F},\widetilde{P})$
such that
\begin{equation}
\label{22.1}
\widetilde{P}^Y\,=\,\P^{Y_t\mid Y_{t-1}=y}=\P^{Y_t\mid Z_{t-1}=z},
\qquad \widetilde{P}^{Y'}\,=\, \P^{Y_t\mid Y_{t-1}=y'}=\P^{Y_t\mid Z_{t-1}=z'},
\end{equation}
and, for
$Z=(X^2,x_1,\ldots,x_{p-1},\lambda,\lambda_1,\ldots,\lambda_{q-1})$
and
$Z'=({X'}^2,x_1',\ldots,x_{p-1}',\lambda',\lambda_1',\ldots,\lambda_{q-1}')$,
\begin{equation}
\label{22.2}
\widetilde{E}\Delta_{\gamma,\delta}(Z,Z') \,\leq\, \kappa \; \Delta_{\gamma,\delta}(z,z')
\end{equation}
holds for some $\kappa <1$.
Actually, according to the model equation (\ref{21.1b}), we have to set
$v=f(x_1^2,\ldots,x_p^2,v_1,\ldots,v_q)$ and
$v'=f({x_1'}^2,\ldots,{x_p'}^2,v_1',\ldots,v_q')$.
Suppose that $(\widetilde{\Omega},\widetilde{\F},\widetilde{P})$ admits the construction of independend
random variables~$U$ and~$V$, both following a uniform distribution on~$[0,1]$.
Let $G_v$ and $G_{v'}$ be the respective distribution functions of $\P^{X_t^2\mid v_t=v}$ and
$\P^{X_t^2\mid v_t=v'}$. We define versions of $X^2$ and ${X'}^2$ by $W:=G_v^{-1}(U)$ and $W':=G_{v'}^{-1}(U)$,
where $G^{-1}$ denotes the generalized inverse
of a generic distribution function~$G$,  $G^{-1}(t)=\inf\{x\colon\; G(x)\geq t\}$.
We still have to determine the signs of $X$ and $X'$, taking into account that the values of
$X^2=W$ and ${X'}^2=W'$ are already determined.
With a view to our proof of absolute regularity, and since the probability $P\big( |X|=|X'|\big)$ is under control,
we will do this in such way that the probability $P\big( |X|=|X'|, X\neq X' \big)$ is as small as possible.
Let $p_l=\P(X_t=l\mid v_t=v)$ and $p_l'=\P(X_t=l\mid v_t=v')$.
If $(p_k\wedge p_k')+(p_{-k}\wedge p_{-k}')\leq P\big(W=W'=k^2\big)$, then we can couple the signs of $X$ and $X'$
such that
\begin{displaymath}
P\big( W=W'=k^2, X=X'=k \big) \,=\, p_k\wedge p_k' \quad \mbox{ and } \quad
P\big( W=W'=k^2, X=X'=-k \big) \,=\, p_{-k}\wedge p_{-k}'.
\end{displaymath}
In this case,
\begin{subequations}
\begin{equation}
\label{22.2a}
P\big( |X|=|X'|=|k|, X\neq X' \big) \,=\, P\big( |X|=|X'|=|k| \big) \,-\, p_k\wedge p_k'
\,-\, p_{-k}\wedge p_{-k}'.
\end{equation}
On the other hand, if $(p_k\wedge p_k')+(p_{-k}\wedge p_{-k}')>P\big(W=W'=k^2\big)$, then we couple
$X$ and $X'$ such that
\begin{displaymath}
P\big( W=W'=k^2, X=X'=\pm k \big) \,=\, p_{\pm k}\wedge p_{\pm k}'
\; \frac{ P\big(W=W'=k^2\big) }{ (p_k\wedge p_k')+(p_{-k}\wedge p_{-k}') },
\end{displaymath}
which leads to
\begin{equation}
\label{22.2b}
P\big( |X|=|X'|=|k|, X\neq X' \big) \,=\, 0.
\end{equation}
\end{subequations}
We denote the corresponding Markov kernels by $\pi^{Y,Y'}$ and $\pi^{Z,Z'}$, respectively.
This contruction produces a pair $(X,X')$ such that
\begin{displaymath}
\widetilde{P}(X=k) \,=\, \P(X_t=k\mid Y_{t-1}=y) \qquad \mbox{ and } \qquad \widetilde{P}(X'=k) \,=\, \P(X_t=k\mid Y_{t-1}=y')
\qquad \forall k\in\Z,
\end{displaymath}
which means that (\ref{22.1}) is satisfied.
Furthermore, it follows from {\bf (A2)} that $X^2$ is stochastically not greater than ${X'}^2$ if $v\leq v$
and vice versa. Since the coupling of these random variables is based on the quantile transform we obtain by (\ref{21.1c})
\begin{equation}
\label{22.3}
\widetilde{E}\big|X^2-{X'}^2\big| \,=\, \big|\widetilde{E}[X^2-{X'}^2]\big| \,=\, |v-v'|.
\end{equation}
Since $P\big(|X|=|X'|=|k|\big)\leq (p_k+p_{-k})\wedge (p_k'\vee p_{-k}')\leq (p_k+p_k'+p_{-k}+p_{-k}')/2$
we obtain from (\ref{22.2a}) and (\ref{22.2b}) that
\begin{eqnarray}
\label{22.4}
P\big( |X|=|X'|, X\neq X' \big)
& = & \sum_{k=1}^\infty P\big( |X|=|X'|=k, X\neq X' \big) \nonumber \\
& \leq & \sum_{k=1}^\infty (p_k+p_k'+p_{-k}+p_{-k}')/2
\,-\, p_k\wedge p_k' \,-\, p_{-k}\wedge p_{-k}' \nonumber \\
& \leq & \sum_{k=-\infty}^\infty \frac{p_k+p_k'}{2} \,-\, p_k\wedge p_k'
\,=\, \frac{1}{2} \, \sum_{k=-\infty}^\infty \big| p_k \,-\, p_k' \big|.
\end{eqnarray}

It follows from (\ref{22.3}) and by {\bf (A1)} that
\begin{eqnarray}
\label{22.11}
\lefteqn{ \widetilde{E} \Delta_{\gamma,\delta}(Z,Z') } \nonumber \\
& = & (\gamma_1+\delta_1) \; \left| f(x_1^2,\ldots,x_p^2,v_1,\ldots,v_q)
\,-\, f({x_1'}^2,\ldots,{x_p'}^2,v_1',\ldots,v_q') \right| \nonumber \\
& & {} \,+\, \sum_{i=2}^p \gamma_i\; |x_{i-1}^2-{x_{i-1}'}^2| + \sum_{j=2}^q \delta_j\; |v_{j-1}-v_{j-1}'| \qquad \\
& \leq & (\gamma_1+\delta_1) \; \left( \sum_{i=1}^p c_i\; |x_i^2-{x_i'}^2| \,+\, \sum_{j=1}^q d_j\; |v_j-v_j'| \right)
\,+\, \sum_{i=2}^p \gamma_i\; |x_{i-1}^2-{x_{i-1}'}^2| \,+\, \sum_{j=2}^q \delta_j\; |v_{j-1}-v_{j-1}'|. \nonumber
\end{eqnarray}
The desired relation of $\widetilde{E}\Delta(Z,Z')\leq\kappa\Delta(z,z')$ would be guaranteed to hold
if we find strictly positive $\gamma_1,\ldots,\gamma_p,\delta_1,\ldots,\delta_p$ such that the
right-hand side of (\ref{22.11}) is less than or equal to
$\kappa\Delta(z,z')=\kappa\big(\sum_{i=1}^p\gamma_i|x_i^2-{x_i'}^2|+\sum_{j=1}^q\delta_j|v_j-v_j'|\big)$,
for all $(x_1^2,\ldots,x_p^2,v_1,\ldots,v_q)$, $({x_1'}^2,\ldots,{x_p'}^2,v_1',\ldots,v_q')\in S^\geq$.
The following lemma provides a bridge from the contraction property {\bf (A1)} for the volatility function
to a contraction property for $Z_t$.
\bigskip
{\lem
\label{L22.1}
Let $c_1,\ldots,c_p,d_1,\ldots,d_q$ be non-negative constants with
\mbox{$\sum_{i=1}^p c_i\,+\, \sum_{j=1}^q d_j<1$}. Then there exist strictly positive constants
$\gamma_1,\ldots,\gamma_p,\lambda_1,\ldots,\lambda_q$ and some $\kappa <1$ such that
\begin{equation}
\label{22.12}
(\gamma_1+\delta_1) \; \Big( \sum_{i=1}^p c_i\; y_i \,+\, \sum_{j=1}^q d_j\; z_j \Big)
\,+\, \sum_{i=2}^p \gamma_i\; y_{i-1} \,+\, \sum_{j=2}^q \delta_j\; z_{j-1}
\,\leq\, \kappa \; \Big( \sum_{i=1}^p \gamma_i\; y_i \,+\, \sum_{j=1}^q \delta_j\; z_j \Big)
\quad
\end{equation}
holds for all $y_1,\ldots,y_p,z_1,\ldots,z_q\geq 0$.
}
\bigskip

Let $\pi^{Y,Y'}$ and $\pi^{Z,Z'}$ be the Markov kernels
which provide the above coupling, that is, for the above pairs of random variables~$(Y,Y')$
and~$(Z,Z')$ we have that $(Y,Y')\sim\pi^{Y,Y'}((y,y'),\cdot)$
and $(Z,Z')\sim\pi^{Z,Z'}((z,z'),\cdot)$, respectively.

The following proposition provides the contraction property which will be instrumental
for the proof of the existence and uniqueness of a stationary distribution as well as for
the derivation of absolute regularity of the count process.
\bigskip
{\prop
\label{P22.1}
Suppose that conditions {\bf (A1)} and {\bf (A2)} are fulfilled.
Let $\gamma_1,\ldots,\gamma_p,\delta_1,\ldots,\delta_p$ and $\kappa<1$ be chosen as in Lemma~\ref{L22.1}. Then
\begin{itemize}
\item[(i)] Let $z,z'\in S^\geq$ be arbitrary. If $(Z,Z')\sim\pi^{Z,Z'}((z,z'),\cdot)$, then
\begin{displaymath}
Z\sim \P^{Z_t\mid Z_{t-1}=z} \qquad \mbox{ and } \qquad Z'\sim \P^{Z_t\mid Z_{t-1}=z'}
\end{displaymath}
and
\begin{displaymath}
\widetilde{E} \Delta_{\gamma,\delta}(Z, Z') \,\leq\, \kappa\; \Delta_{\gamma,\delta}(z, z').
\end{displaymath}
\item[(ii)] Let $((\widetilde{Z}_t,\widetilde{Z}_t'))_{t\in\Z}$ be a Markov chain on $(\widetilde{\Omega},\widetilde{\F},\widetilde{P})$
with transition kernel~$\pi^{Z,Z'}$. Then
\bd
\widetilde{E} \Delta_{\gamma,\delta}\left(\widetilde{Z}_t,\widetilde{Z}_t'\right)
\,\leq\, \kappa \; \widetilde{E} \Delta_{\gamma,\delta}\left(\widetilde{Z}_{t-1},\widetilde{Z}_{t-1}'\right).
\ed
\end{itemize}
}
\bigskip

In order to derive stationarity properties of the process $(Z_t)_{t\in\Z}$,
we further translate the contraction result in Proposition~\ref{P22.1} into
a contraction property of the corresponding distributions.
For the metric $\Delta_{\gamma,\delta}$ on~$S^\geq$, we define
\begin{displaymath}
{\mathcal P}(S^\geq) \,=\, \big\{Q\colon\quad Q \mbox{ is a probability distribution on } S^\geq
\mbox{ with } \int \Delta_{\gamma,\delta}(z_0,z)\, dQ(z)<\infty \big\},
\end{displaymath}
where $z_0\in S^\geq$ is arbitrary.
For two probability measures~$Q,Q'\in {\mathcal P}(S^\geq)$, we define the Kantorovich distance
based on the metric $\Delta_{\gamma,\delta}$ (also known as Wasserstein $L^1$ distance)
by
\bd
{\mathcal K}(Q, Q') \,:=\, \inf_{Z\sim Q, Z'\sim Q'} \widetilde{E} \Delta(Z, Z'),
\ed
where the infimum is taken over all random variables~$Z$ and~$Z'$ defined on a common
probability space $(\widetilde{\Omega},\widetilde{\F},\widetilde{P})$ with respective laws~$Q$ and~$Q'$.
We denote the Markov kernel of the processes $(Y_t)_{t\in\Z}$ and $(Z_t)_{t\in\Z}$ by $\pi^Y$ and $\pi^Z$,
respectively.
The following result follows immediately from Proposition~\ref{P22.1}.
\bigskip
{\prop
\label{P22.2}
Suppose that conditions {\bf (A1)} and {\bf (A2)} are fulfilled.
Let $Q,Q'\in {\mathcal P}(S^\geq)$ be arbitrary distributions. Then, for $\kappa<1$
given in Lemma~\ref{L22.1},
\bd
{\mathcal K}(Q \pi_\theta^Z, Q' \pi_\theta^Z) \,\leq\, \kappa \, {\mathcal K}(Q, Q').
\ed
}
\bigskip

%%%%%%%%%%%%%%%%%%%%%%%%%%%%%%%%%%%%%%%%%%%%%%%%%%%%%%%%%%%%%%%%%%%%%%%%%%%%%%%
\subsection{Existence and uniqueness of a stationary distribution}
\label{SS2.3}
%%%%%%%%%%%%%%%%%%%%%%%%%%%%%%%%%%%%%%%%%%%%%%%%%%%%%%%%%%%%%%%%%%%%%%%%%%%%%%%

Proposition~\ref{P22.2} shows that the mapping $\pi^Z$ is contractive. Therefore, we can
conclude by the Banach fixed point theorem that the Markov process $(Z_t)_{t\in\Z}$
has a unique stationary distribution. A simple extra argument shows that that the
Markov process $(Y_t)_{t\in\Z}$ has this property as well.

\bigskip
{\thm
\label{T23.1}
Suppose that conditions {\bf (A1)} and {\bf (A2)} are fulfilled.
\begin{itemize}
\item[(i)]
The Markov process $(Z_t)_{t\in\Z}$ with transition kernel $\pi^Z$ has a unique stationary distribution $Q^Z$.
For $Z_0=(X_0^2,\ldots,X_{1-p}^2,v_0,\ldots,v_{1-q})$, we have that
\begin{equation}
\label{23.1}
E\left[ X_0^2 \,+\, v_0 \right] \,<\, \infty.
\end{equation}
\item[(ii)]
The Markov process $(Y_t)_{t\in\Z}$ with transition kernel $\pi^Y$ has a unique stationary distribution $Q^Y$.
\end{itemize}
}
\bigskip

\begin{rem}
\label{R2.1}
The reader might wonder why we don't derive weak dependence properties introduced by \citet{DL99}. Indeed e.g.
\citet{DN08} describe statistical procedures where mixing can be replaced by weak dependence conditions.
If $X\sim Q_v=\Skellam(v/2,v/2)$ is symmetric, then $EX=0$ and $EX^2=v$.
Therefore it is natural to model the volatility process as in (\ref{21.1b}),
where the volatilities appear linearly while the count variables are squared.
The properties of Skellam models make also natural the inhomogeneity of \eqref{21.1a} and \eqref{21.1b}
which include both linear and squared factors.
In the simplest  case of a SkellamARCH(1)-process, with $v_t  =  f(X_{t-1}^2)$,
the function $x\mapsto f(x^2)$ may not be Lipschitz and thus contraction does not hold.
Anyway, the process $Y_t=X_t^2$ is again contractive if  $\mbox{Lip }f<1$.
Symmetry of the distribution  $Q_v$ implies that   $X_t=\sigma_t\sqrt{Y_t}$ is a solution of  \eqref{21.1a} and \eqref{21.1b}   if
$(\sigma_t)_t$ is an iid sequence of symmetric signs ($P(\sigma_t=\pm1)=1/2$).
Then $\tau-$dependence of the process $(Y_t)$ follows as in \citet{DW08}.
Now, since  the 1-Lipschitz function $g(y)=y\wedge \sqrt{y}$   equals $y\mapsto \sqrt{y}$ on $\N_0$, then $X_t=\sigma_tg(Y_t)$;
heredity properties of weak dependence imply   geometric $\tau-$dependence of  $(X_t)_t$; see \citet{DDLLLP07}.
We proved $\tau-$dependence in this very special symmetric case; in order to work in  a more general setting we  switch
in Subsection~\ref{SS2.4}  to the more  standard  $\beta-$mixing condition to derive asymptotic theory for the statistical analysis.
\end{rem}

%%%%%%%%%%%%%%%%%%%%%%%%%%%%%%%%%%%%%%%%%%%%%%%%%%%%%%%%%%%%%%%%%%%%%%%%%%%%%%%
\subsection{Absolute regularity}
\label{SS2.4}
%%%%%%%%%%%%%%%%%%%%%%%%%%%%%%%%%%%%%%%%%%%%%%%%%%%%%%%%%%%%%%%%%%%%%%%%%%%%%%%

For the related case of Poisson count processes with a GARCH-type structure,
absolute regularity has been first proved for contractive INGARCH(1,1) processes
in \citet{Neu11}.
This has been generalized in \citet{DN19} to semi-contractive models and
in \citet{DLN20} to the case of possibly non-stationary processes.
In all of these papers, the mixing properties were derived by an explicit
coupling of two versions of the processes which were tailor-made for the respective
properties of the processes.
In the current work, our approach is slightly different. We derive both stationarity
and mixing properties on the basis of a one-step contractivity property
given in Proposition~\ref{P22.1}.

Let $(\Omega,{\mathcal A},P)$ be a probability space and ${\mathcal A}_1$, ${\mathcal A}_2$
be two sub-$\sigma$-algebras of ${\mathcal A}$. Then the coefficient of absolute regularity is defined as
\bd
\beta({\mathcal A}_1,{\mathcal A}_2)
\,=\, E\left[ \sup\left\{ |P(B\mid {\mathcal A}_1) \,-\, P(B)|\colon \;\; B\in {\mathcal A}_2 \right\} \right].
\ed
For a strictly stationary process ${\mathbf Y}=(Y_t)_t$ on $(\Omega,{\mathcal A},P)$, the
coefficients of absolute regularity are defined as
\bd
\beta^Y(n) \,=\, \beta\left( \sigma(Y_0,Y_{-1},\ldots), \sigma(Y_n,Y_{n+1},\ldots) \right).
\ed
For the count process $(X_t)_t$ on $(\Omega,\F,\P)$, we obtain the following estimate
of the coefficients of absolute regularity.
\bea
\label{24.1}
\lefteqn{ \beta^X(n) } \nonumber \\
& = & \beta\big( \sigma(X_0,X_{-1},\ldots), \sigma(X_n,X_{n+1},\ldots) \big) \nonumber \\
& \leq & \beta\big( \F_0, \sigma(X_n,X_{n+1},\ldots) \big) \nonumber \\
& = & \beta\big( \sigma(Y_0), \sigma(X_n,X_{n+1},\ldots) \big) \nonumber \\
& = & \E\left[ \sup_{C\in\sigma({\mathcal Z})} \left\{ \Big| \P_\theta\big((X_n,X_{n+1},\ldots)\in C\mid Z_k\big)
\,-\, \P_\theta\big((X_n,X_{n+1},\ldots)\in C\big) \Big| \right\} \right], \qquad \quad
\eea
where ${\mathcal Z}=\{A\times\Z\times\Z\times\cdots\mid A\subseteq\Z^m, m\in\N\}$
is the system of cylinder sets.
At this point we employ a coupling argument.
Let $((\widetilde{Y}_t,\widetilde{Y}_t'))_{t\in\N_0}$ be a Markov chain on a probability space
$(\widetilde{\Omega},\widetilde{\F},\widetilde{P})$ with transition kernel $\pi^{Y,Y'}$
and independent variables $\widetilde{Y}_0,\widetilde{Y}_0'\sim\P_\theta^{Z_k}$.
Then
\bea
\label{24.2}
\lefteqn{ \E\left[ \sup_{C\in\sigma({\mathcal Z})} \left\{ \Big| \P_\theta\big((X_k,X_{n+1},\ldots)\in C\mid Z_k\big)
\,-\, \P_\theta\big((X_k,X_{n+1},\ldots)\in C\big) \Big| \right\} \right] } \nonumber \\
& \leq & \widetilde{E}\Big[ \sup_{C\in\sigma({\mathcal Z})} \big\{ \big|
\widetilde{P}\left( (\widetilde{X}_n,\widetilde{X}_{n+1},\ldots)\in C\mid \widetilde{Y}_0 \right)
\,-\, \widetilde{P}\left( (\widetilde{X}_n',\widetilde{X}_{n+1}',\ldots)\in C\mid \widetilde{Y}_0' \right) \big| \big\} \Big]
\nonumber \\
& \leq & \widetilde{P}\left( \widetilde{X}_{n+k} \neq \widetilde{X}_{n+k}' \quad \mbox{ for some } k\geq 0 \right)
\qquad \qquad \qquad \qquad \qquad \qquad
\nonumber \\
& \leq & \sum_{k=0}^\infty \widetilde{P}\left( \widetilde{X}_{n+k} \neq \widetilde{X}_{n+k}' \right).
\eea

At this point we will more closely examine the remaining  part of our approach to derive upper estimates for the mixing coefficients.
If the count variables are non-negative, then $\widetilde{X}_{n+k}=\widetilde{X}_{n+k}'$ is equivalent to
$|\widetilde{X}_{n+k}|=|\widetilde{X}_{n+k}'|$.
Moreover, if the probability mass functions of the~$Q_v$ are symmetric about zero, then (\ref{22.2a}) and (\ref{22.2b}) ensure that
$\widetilde{X}_{n+k}$ and $\widetilde{X}_{n+k}'$ have always the same sign which means again that
$\widetilde{X}_{n+k}=\widetilde{X}_{n+k}'$ is equivalent to $|\widetilde{X}_{n+k}|=|\widetilde{X}_{n+k}'|$.
In both cases, we conclude from (\ref{24.1}) and (\ref{24.2}) that
\begin{eqnarray}
\label{24.3}
\beta^X(n) & \leq & \sum_{k=0}^\infty \widetilde{P}\left( \widetilde{X}_{n+k} \neq \widetilde{X}_{n+k}' \right) \nonumber \\
& \leq & \frac{1}{\gamma_1} \; \sum_{k=0}^\infty \widetilde{E} \Delta_{\gamma,\delta}(\widetilde{Z}_{n+k}, \widetilde{Z}_{n+k}') \nonumber \\
& \leq & \frac{1}{\gamma_1} \; \frac{\kappa^n}{1-\kappa} \; \widetilde{E} \Delta_{\gamma,\delta}(\widetilde{Z}_0, \widetilde{Z}_0').
\end{eqnarray}
Otherwise, we assume that $(Q_v)_{v\in V}$ is such that, for some $K<\infty$,
\begin{equation}
\label{24.4}
\frac{1}{2}\, \sum_{k=-\infty}^\infty \big| Q_v(\{k\}) \,-\, Q_{v'}(\{k\}) \big| \,\leq\, K\,|v-v'| \qquad \forall v,v'\in V.
\end{equation}
Then we obtain by (\ref{22.4}) that
\begin{eqnarray*}
\widetilde{P}\big( \widetilde{X}_{n+k} \neq \widetilde{X}_{n+k}' \big)
& = & \widetilde{P}\big( |\widetilde{X}_{n+k}| \neq |\widetilde{X}_{n+k}'| \big)
\,+\, \widetilde{P}\big( |\widetilde{X}_{n+k}| = |\widetilde{X}_{n+k}'|, \widetilde{X}_{n+k} \neq \widetilde{X}_{n+k}' \big) \\
& \leq & \widetilde{P}\big( |\widetilde{X}_{n+k}| \neq |\widetilde{X}_{n+k}'| \big)
\,+\, K\; \big| v_{n+k} \,-\, v_{n+k}' \big|.
\end{eqnarray*}
In this case, we obtain that
\begin{eqnarray}
\label{24.5}
\beta^X(n) & \leq & \left( \frac{1}{\gamma_1} \,+\, \frac{K}{\delta_1} \right)
\sum_{k=0}^\infty \Delta_{\gamma,\delta}(\widetilde{Z}_{n+k}, \widetilde{Z}_{n+k}') \nonumber \\
& \leq & \left( \frac{1}{\gamma_1} \,+\, \frac{K}{\delta_1} \right) \; \frac{\kappa^n}{1-\kappa}
\; \widetilde{E} \Delta_{\gamma,\delta}(\widetilde{Z}_0, \widetilde{Z}_0').
\end{eqnarray}
\bigskip

{\thm
\label{T24.1}
Suppose that conditions {\bf (A1)} and {\bf (A2)} are fulfilled and that the process
$(Y_t)_{t\in\Z}$ is stationary.
Furthermore we assume that $(Q_v)_{v\in V}$ satisfies one of the following conditions.
\begin{itemize}
\item[a)\quad] $Q_v(\N_0)=1\quad \forall v\in V$,
\item[b)\quad] the probability mass functions of $(Q_v)_{v\in V}$ are symmetric about zero,
\item[c)\quad] (\ref{24.4}) is fulfilled for some $K<\infty$.
\end{itemize}
Then there exists some $\rho<1$ such that
\bd
\beta^X(n) \,=\, O\left( \rho^n \right).
\ed
}
\begin{rem}
\label{R2.2}
The results of our paper are heavily based on the (fully) contractive
condition {\bf (A1)}
on the volatility function~$f$.
In a related work, \citet{DN19}, a weaker so-called semi-contractive
condition,
\begin{displaymath}
\big| f(x_1,\ldots,x_p,\lambda_1,\ldots,\lambda_q) \,-\,
f(x_1,\ldots,x_p,\lambda_1',\ldots,\lambda_q') \big|
\,\leq\, \sum_{j=1}^q d_j |\lambda_j - \lambda'|,
\end{displaymath}
for some non-negative $d_1,\ldots,d_q$ such that $\sum_{j=1}^q d_j<1$,
was imposed which then resulted a a slower subexponential decay of the
coefficients of
absolute regularity.
In our context, it seems also be possible to derive properties such as
existence and uniqueness
of a stationary distribution and absolute regularity under a
semi-contractive condition
if some appropriate drift condition is added.
Without any kind of contractivity condition, the approach used in this
paper fails and there are
counterexamples showing that then our results are non longer valid.
Consider the special case of the linear model (\ref{LGARCH}), where
$\omega>0$ and
$\alpha_1,\ldots,\alpha_p,\beta_1,\ldots,\beta_q$ are non-negative with
$\sum_{i=1}^p\alpha_i+\sum_{j=1}^q\beta_j\geq 1$.
Then
\begin{displaymath}
\E v_t \,=\, \omega \,+\, \sum_{i=1}^p \alpha_i \E v_{t-i} \,+\,
\sum_{j=1}^q \beta_j \E v_{t-j},
\end{displaymath}
which shows that a stationary process $((X_t,v_t))_{t\in\Z}$
satisfying (\ref{21.1a}) to (\ref{21.1c})
does not exist.
\end{rem}
\bigskip

%%%%%%%%%%%%%%%%%%%%%%%%%%%%%%%%%%%%%%%%%%%%%%%%%%%%%%%%%%%%%%%%%%%%%%%%%%%%%%
\section{Applications}
\label{S3}
%%%%%%%%%%%%%%%%%%%%%%%%%%%%%%%%%%%%%%%%%%%%%%%%%%%%%%%%%%%%%%%%%%%%%%%%%%%%%%%
We choose to develop the asymptotic theory for the OLSE of Skellam models \eqref{3.0}
as the most standard application of the above results. Much more may be done including tests of goodness-of-fit
as in \cite{DLN20}. Prediction or model selection issues are also important and should be developed
theoretically.
Additional research work will make use of the bound of absolute regularity for many other
questions such a more quantitative study of prediction, qualitative tests of goodness-of-fit
such as model choice problems, or more nonparametric based statistics or resampling or subsampling procedures.
%%%%%%%%%%%%%%%%%%%%%%%%%%%%%%%%%%%%%%%%%%%%%%%%%%%%%%%%%%%%%%%%%%%%%%%%%%%%%%%
\subsection{OLSE of a Skellam-ARCH model}
\label{SS3.1}
%%%%%%%%%%%%%%%%%%%%%%%%%%%%%%%%%%%%%%%%%%%%%%%%%%%%%%%%%%%%%%%%%%%%%%%%%%%%%%%

We consider the special case of an Skellam-ARCH($p$) model, where (\ref{21.1b})
reduces to
\begin{equation}
\label{3.0}
v_t \,=\, \omega \,+\, \sum_{i=1}^p \alpha_i \; X_{t-i}^2.
\end{equation}
We assume that $\omega>0$, and that $\alpha_1,\ldots,\alpha_p$ are non-negative
with $\alpha=\sum_{i=1}^p\alpha_i<1/\sqrt{3}$.
We further assume that the process $((X_t,v_t))_{t\in\Z}$ is in its
unique stationary regime. On the basis of observations $X_1,\ldots,X_n$, we
intend to estimate the vector of unknown parameters $\theta=(\omega,\alpha_1,\ldots,\alpha_p)^T$.
We embed the observed random variables into a linear regression model,
\begin{displaymath}
X_t^2 \,=\, \omega \,+\, \sum_{i=1}^p X_{t-i}^2 \; \alpha_i \,+\, \varepsilon_t,
\qquad t=p+1,\ldots,n,
\end{displaymath}
where $\varepsilon_t=X_t^2-v_t$ satisfies $E(\varepsilon_t\mid \F_{t-1})=0$ a.s.
Then the ordinary least squares estimator is given by
\begin{eqnarray*}
\widehat{\theta}_n
& \in & \arg \min_\theta \sum_{t=p+1}^n
\left( X_t^2 \;\; - \;\; ( \omega \,+\, \alpha_1\; X_{t-1}^2 \,+\, \cdots \,+\, \alpha_p\; X_{t-p}^2) \right)^2 \\
& = & \arg \min_\theta \left\| Y_{(n)} \,-\, X_{(n)} \theta \right\|^2,
\end{eqnarray*}
where
\begin{displaymath}
Y_{(n)} \,=\, \left( \begin{array}{c}
X_{p+1}^2 \\ \vdots \\ X_n^2
\end{array} \right), \qquad
X_{(n)} \,=\, \left( \begin{array}{cccc}
1 & X_p^2 & \hdots & X_1^2 \\
\vdots & \vdots & \ddots & \vdots \\
1 & X_{n-1}^2 & \hdots & X_{n-p}^2
\end{array} \right).
\end{displaymath}
If the matrix $X_{(n)}^T X_{(n)}$ is regular, then $\widehat{\theta}_n$ is uniquely defined and
\begin{equation}\label{est}
\widehat{\theta}_n \,=\, \left( X_{(n)}^T X_{(n)} \right)^{-1} \;
X_{(n)}^T Y_{(n)},
\end{equation}
which implies that
\begin{equation}
\label{3.1}
\sqrt{n} \; \left( \widehat{\theta}_n \,-\, \theta \right)
\,=\, \left( \frac{1}{n} X_{(n)}^T X_{(n)} \right)^{-1} \; \frac{1}{\sqrt{n}} \; X_{(n)}^T \varepsilon_{(n)},
\end{equation}
where $\varepsilon_{(n)}=(\varepsilon_{p+1},\ldots,\varepsilon_n)^T$.

The condition $\sum_{i=1}^p\alpha_i<1/\sqrt{3}$ ensures by Lemma~\ref{L4.1} that $E X_0^4<\infty$.
Hence, we obtain from the ergodic theorem that
\begin{equation}
\label{3.2}
\frac{1}{n} X_{(n)}^T X_{(n)} \,\stackrel{a.s.}{\longrightarrow}\, \Sigma
\,=\, \left( \begin{array}{cccc}
1 & E X_1^2 & \hdots & E X_p^2 \\
E X_1^2 & E X_1^2 X_1^2 & \hdots & E X_1^2 X_p^2 \\
\vdots & \vdots & \ddots & \vdots \\
E X_p^2 & E X_p^2 X_1^2 & \hdots & E X_p^2 X_p^2
\end{array} \right).
\end{equation}
Lemma~\ref{L4.2} below shows that $\Sigma$ is a regular matrix which means that equation (\ref{3.1})
holds true with a probability tending to~1.
Furthermore, it follows from a central limit theorem for sums of martingale differences (Corollary 3.1 in \cite{HH80}  page 58) and the Cram\'er-Wold device
that
\begin{equation}
\label{3.3}
W_n \,:=\, \frac{1}{\sqrt{n}} \; X_{(n)}^T \varepsilon_{(n)} \,\stackrel{d}{\longrightarrow}\, Z_0,
\end{equation}
where $Z_0 \sim {\mathcal N}( 0_{p+1}, \eta^2\;\Sigma )$ and $\eta^2=E\varepsilon_t^2$.

{}From (\ref{3.1}) to (\ref{3.3}) we conclude that
\begin{equation}\label{CLT_MLE}
\sqrt{n} \; \left( \widehat{\theta}_n \,-\, \theta \right) \,\stackrel{d}{\longrightarrow}\, Z
\,\sim\, {\mathcal N}( 0_{p+1}, \eta^2\; \Sigma^{-1} ).
\end{equation}
%\medskip

%%%%%%%%%%%%%%%%%%%%%%%%%%%%%%%%%%%%%%%%%%%%%%%%%%%%%%%%%%%%%%%%%%%%%%%%%%%%%%%
\subsection{Simulation study}
\label{SS3.2}
%%%%%%%%%%%%%%%%%%%%%%%%%%%%%%%%%%%%%%%%%%%%%%%%%%%%%%%%%%%%%%%%%%%%%%%%

We simulate a process $((X_t,v_t))_t$, where $v_t$ obeys (\ref{3.0}) and
$X_t\mid \F_t\sim Q_{v_t}=\Skellam(v_t/2,v_t/2)$. The parameters in (\ref{3.0}) are chosen such that
$\omega>0$ and $\sum\alpha_i<1/\sqrt{3}$, which ensures finiteness of fourth moments of the count variables;
see Lemma~\ref{L4.1} below. We assume a suitable set of values for the different order $p$,
$\alpha_{1}=0.26$, $\alpha_{2}=0.16$, $\alpha_{3}=0.11$, $\alpha_{4}=0.02$,
for sample sizes $T=30,80,100,500$, and 1000.
1000 replication are made for each sample size and the simulated mean estimates and their
corresponding standard errors as deduced from the result~\eqref{CLT_MLE}.

\begin{table}[H]
{\small
\begin{tabular}{|c|c|c|c|c|c|c|}
\hline
$p$	&   $T$	& $\omega=1.50$	& $\alpha_{1}=0.26$&	$\alpha_{2}=0.16$	& $\alpha_{3}=0.11$ &  $ \alpha_{4}=0.02$ \\ \hline
1	&30	 &  1.251  (0.321)&	0.178  (0.211)&	&	&     \\ \hline
	&80	 &1.355  (0.151)	& 0.225 (0.188)&	&	&     \\  \hline
	&100	&1.751 (0.101)&	0.231 (0.124)&		&  &\\  \hline
	&500	&1.442	(0.087)& 0.244 	(0.091)	&	&   & \\ \hline
	&1000	&1.542 (0.075)&	0.136 (0.088)&		& & \\ \hline
\hline		
2	& 30	&1.389 (0.278) &	0.232 (0.209)&	0.152 	(0.327)& &	    \\ \hline
	&80	&1.477 (0.150)&	0.246 (0.123)	&0.166 (0.111)&	 &  \\ \hline
	&100	&1.511 (0.081)	&0.271 (0.099)&	0.148 (0.098)	& & \\ \hline
	&500	&1.552 (0.032)	&0.276 (0.042)&	0.152 (0.038)	& &  \\ \hline
	&1000	&1.467 (0.022)&	0.255 (0.031)&	0.169 	(0.021)	& &\\ \hline
 \hline					
3	& 30 &  1.481   (0.455)  & 0.244 (0.303)	  & 0.152 (0.276)  &	0.104  (0.152)  &    \\ \hline
	&80	&  1.551 (0.210)  &	0.232 (0.155)	& 0.166 (0.101)  &	0.119 (0.110) & \\ \hline
	&100&  1.462   (0.111) & 0.255 (0.101)   & 0.147 (0.089)	 & 0.114   (0.088)& \\ \hline
	&500&  1. 541 (0.088) & 0.275 (0.076)   & 0.158 (0.042)	& 0.121 (0.034) &  \\ \hline
	&1000	&1.4601 (0.061)&	0.266 (0.045) &	0.166	(0.034) &0.111 (0.026) &  \\ \hline
 \hline					
4	& 30 &  1.495 (0.323) & 0.255 (0.212)	  & 0.152 (0.318) &	0.112 (0.176)  & 0.018  (0.272)  \\ \hline
	&80	&  1.510  (0.188)&	0.276 (0.124)	& 0.164 (0.232)  &	0.114  (0.123)&0.015  (0.103)\\ \hline
	&100&  1.498  (0.092)& 0.237 (0.110)   & 0.166	  (0.101)& 0.110 (0.075)  &0.019  (0.064)  \\ \hline
	&500&  1.502  (0.088)& 0.242  (0.064) & 0.167	 (0.054)& 0.113  (0.033) &0.017   (0.042) \\ \hline
	&1000	&1.489 (0.052)&	0.262 (0.043) &	0.159 (0.038)	&0.109  (0.018)& 0.018  (0.028)	 \\ \hline

\end{tabular}
}
\caption{Simulated Mean estimates and standard errors/100 Simulations}\label{table1}
\end{table}
\noindent
The estimates reflect that for increased sample size, the values of the different parameters
become more consistent with the standard errors that are seen to be constantly decreasing; see Table \ref{table1}.
It is worth reporting that for some simulation processes, the standard OLS equation (\ref{3.1})
does not yield relevant output estimates since the constraints on the $\widehat{\alpha}_{i}$ and $\widehat{\omega}$
were not conformed. To overcome this shortcomings, we apply the {\tt QP.solve} routine with the appropriate
constraint matrix to obtain reliable results. Furthermore, for $p=4$, some simulations initially failed since
the Hessian matrix was near to ill-conditioned and the evaluation of the inverse was then suitably handled
by the {\tt ginv} function. Broadly, after accommodating these computational amendments,
the average number of convergent simulations turn around 92 \%, 87 \%, 85 \% and 75 \% for $p=1,2,3,4$ respectively.
As noticed from these percentages, as we increase the order $p$, we expect some number of failed simulations.

%%%%%%%%%%%%%%%%%%%%%%%%%%%%%%%%%%%%%%%%%%%%%%%%%%%%%%%%%%%%%%%%%%%%%%%%%%%%%%%%%%%%%%%%%
\section{Proofs and some auxiliary results}
\label{S4}
\subsection{Proofs of the main results}
\label{SS4.1}
%%%%%%%%%%%%%%%%%%%%%%%%%%%%%%%%%%%%%%%%%%%%%%%%%%%%%%%%%%%%%%%%%%%%%%%%%%%%%%%%%%%%%%%%%

\begin{proof}[Proof of Lemma~\ref{L22.1}]
A comparison of coefficients in (\ref{22.12}) reveals that it suffices to find strictly positive constants
$\gamma_1,\ldots,\gamma_p,\delta_1,\ldots,\delta_q$ such that the following inequalities are satisfied.
\bea
\label{pl221.1}
\big(\gamma_1 \,+\, \delta_1\big) \; c_1 \,+\, \gamma_2 & \leq & \kappa\; \gamma_1 \nonumber \\
& \vdots & \nonumber \\
\big(\gamma_1 \,+\, \delta_1\big) \; c_{p-1} \,+\, \gamma_p & \leq & \kappa\; \gamma_{p-1} \nonumber \\
\big(\gamma_1 \,+\, \delta_1\big) \; c_p & \leq & \kappa\; \gamma_p \nonumber \\
\big(\gamma_1 \,+\, \delta_1\big) \; d_1 \,+\, \delta_2 & \leq & \kappa\; \delta_1 \nonumber \\
& \vdots & \nonumber \\
\big(\gamma_1 \,+\, \delta_1\big) \; d_{q-1} \,+\, \delta_q & \leq & \kappa\; \delta_{q-1} \nonumber \\
\big(\gamma_1 \,+\, \delta_1\big) \; d_q & \leq & \kappa\; \delta_q.
\eea
We set, w.l.o.g., $\gamma_1+\delta_1=1$. Let $\epsilon=(1-L)/(p+q)$.
We consider the following system of equations.
\bean
c_p \,+\, \epsilon & = & \gamma_p \\
c_{p-1} \,+\, \gamma_p \,+\, \epsilon & = & \gamma_{p-1} \\
& \vdots & \\
c_1 \,+\, \gamma_2 \,+\, \epsilon & = & \gamma_1 \\
d_q \,+\, \epsilon & = & \delta_q \\
d_{q-1} \,+\, \delta_q \,+\, \epsilon & = & \delta_{q-1} \\
& \vdots & \\
d_1 \,+\, \delta_2 \,+\, \epsilon & = & \delta_1.
\eean
It is obvious that this system of equations has a unique solution
with strictly positive $\gamma_1,\ldots,\gamma_p,\delta_1,\ldots,\delta_q$.
Moreover, it follows from
\bd
\sum_{i=1}^p c_i \,+\, \sum_{j=1}^q d_j \,+\, \sum_{i=2}^p \gamma_i \,+\, \sum_{j=2}^q \delta_j
\,+\, (p+q)\; \epsilon \,=\, \sum_{i=1}^p \gamma_i \,+\, \sum_{j=1}^q \delta_j
\ed
that $\gamma_1+\delta_1=1$, as required.
Therefore, we see that, with such a choice of $\gamma_1,\ldots,\gamma_p,\delta_1,\ldots,\delta_q$,
the following strict inequalities are fulfilled.
\bean
c_1 \,+\, \gamma_2 & < & \gamma_1 \\
& \vdots & \\
c_{p-1} \,+\, \gamma_p & < & \gamma_{p-1} \\
c_p & < & \gamma_p \\
d_1 \,+\, \delta_2 & < & \delta_1 \\
& \vdots & \\
d_{q-1} \,+\, \delta_q & < & \delta_{q-1} \\
d_q & < & \delta_q. \\
\eean
Choosing $\kappa=\max\{(c_1+\gamma_2)/\gamma_1,\ldots,(c_{p-1}+\gamma_p)/\gamma_{p-1},c_p/\gamma_p,
(d_1+\delta_2)/\delta_1,\ldots,(d_{q-1}+\delta_q)/\delta_{q-1},d_q/\delta_q\}$ we obtain that
the system of inequalities (\ref{pl221.1}) is satisfied.
\end{proof}
\bigskip

\begin{proof}[Proof of Proposition~\ref{P22.1}]
(i) follows from (\ref{22.11}) and (\ref{22.12}), and
(ii) is an immediate consequence of (i).
\end{proof}
\bigskip

\begin{proof}[Proof of Proposition~\ref{P22.2}]
Let $Q$ and $Q'$ be arbitrary probability measures supported in~$S^\geq$ and let~$\xi$
be the optimal coupling of~$Q$ and~$Q'$ w.r.t.~the Kantorovich distance, that is,
\bd
{\mathcal K}(Q,Q')
\,=\, \int_{S^\geq\times S^\geq} \Delta_{\gamma,\delta}(z,z') \, \xi(dz, dz').
\ed
Then $\xi\pi^{Z,Z'}$ is a coupling of $Q\pi^Z$ and $Q'\pi^Z$ and it follows
from Proposition~\ref{P2.1}(i) that
\bean
{\mathcal K}(Q\pi, Q'\pi ) & \leq & \int \Delta_{\gamma,\delta}(u,u') \, \xi\pi^{Z,Z'}(du, du') \\
& = & \int \left[ \int \Delta_{\gamma,\delta}(u,u') \, \pi^{Z,Z'}((z,z'), du\, du') \right]
\, \xi(dz, dz') \\
& \leq & \kappa \; \int \Delta_{\gamma,\delta}(z,z') \, \xi(dz, dz') \,=\, \kappa \; {\mathcal K}(Q,Q').
\eean
\end{proof}
\medskip

\begin{proof}[Proof of Theorem~\ref{T23.1}]
We consider first the Markov process $(Z_t)_{t\in\Z}$.
Let
\bd
{\mathcal P} \,=\, \left\{Q\colon\;\; Q \mbox{ is a probability distribution based in } S^\geq,
\int_S\;\displaystyle  \sum_{i=1}^{2p}\;|x_i|\, Q(dx)<\infty \right\}.
\ed
It is well known that the space ${\mathcal P}$ equipped with the Kantorovich metric~${\mathcal K}$ is complete.
Since by Proposition~\ref{P22.2} the mapping $\pi^Z$ is contractive it follows by the Banach fixed point
theorem that the Markov kernel $\pi^Z$ admits a unique fixed point~$Q^Z$, i.e. $Q^Z\pi^Z=Q^Z$.
In other words, $Q^Z$ is the unique stationary distribution of the process $(Z_t)_{t\in\Z}$.

Now we consider the process $(Y_t)_{t\in\Z}$.
If the $X_t$ are non-negative random variables, then we have a one-to-one relationship
between~$Z_t$ and~$Y_t$ and, for $Z_0=\big(X_0^2,\ldots,X_{1-p}^2,v_0,\ldots,v_{1-q}\big)\sim Q^Z$,
the distribution of the vector $Y_0=\big(X_0,\ldots,X_{1-p},v_0,\ldots,v_{1-q}\big)$
is the unique stationary distribution of $(Y_t)_{t\in\Z}$.

If $X_t$ attains both positive and negative values, we need a simple extra argument.
Suppose that $Z_0=\big(X_0^2,\ldots,X_{1-p}^2,v_0,\ldots,v_{1-q}\big)\sim Q^Z$.
Now we can recursively generate suitable $(v_1,X_1),\ldots,(v_p,X_p)$ as follows.
We set $v_1=f\big(X_0^2,\ldots,X_{1-p}^2,v_0,\ldots,v_{1-q}\big)$ and generate $X_1\sim Q_{v_1}$.
Then we set $v_2=f\big(X_1^2,\ldots,X_{-p}^2,v_1,\ldots,v_{-q}\big)$ and choose
$X_1\sim Q_{v_1}$, and so on.
After~$p$ such steps we have collected enough $X_t$s with suitable signs and the
random vector $Y_p:=\big(X_p,\ldots,X_1,v_p,\ldots,v_{p-q+1}\big)$ has the unique
stationary distribution, say $Q^Y$, of $(Y_t)_{t\in\Z}$.
\end{proof}
\bigskip

\begin{proof}[Proof of Theorem~\ref{T24.1}]
Since by (\ref{23.1}) $\widetilde{E} \Delta_{\gamma,\delta}(\widetilde{Z}_0,\widetilde{Z}_0')<\infty$,
this theorem is an immediate consequence of (\ref{24.3}) and (\ref{24.5}).
\end{proof}

%%%%%%%%%%%%%%%%%%%%%%%%%%%%%%%%%%%%%%%%%%%%%%%%%%%%%%%%%%%%%%%%%%%%%%%%%%%%%%%%%%%%%%%%%
\subsection{Some auxiliary results}
\label{SS4.2}
%%%%%%%%%%%%%%%%%%%%%%%%%%%%%%%%%%%%%%%%%%%%%%%%%%%%%%%%%%%%%%%%%%%%%%%%%%%%%%%%%%%%%%%%%

{\lem
\label{L4.1}
Let $((X_t,v_t))_{t\in\Z}$ be a stationary process satisfying (\ref{3.0}),
where $\omega,\alpha_1,\ldots,\alpha_p$ are non-negative constants
and let $X_t\mid \F_{t-1}\sim\Skellam(v_t/2,v_t/2)$.
\begin{itemize}
\item[(i)] If $\displaystyle \alpha=\sum_{i=1}^n \alpha_i<1$, then
$E X_0^2 \,<\, \infty$.
\item[(ii)] If $\displaystyle \alpha=\sum_{i=1}^n \alpha_i<\frac{1}{\sqrt{3}}$, then
$E X_0^4 \,<\, \infty$.
\end{itemize}
}
\bigskip

\begin{proof}[Proof of Lemma~\ref{L4.1}]
Let $((\widetilde{X}_t,\widetilde{v}_t))_{t\in\N}$ be Skellam-ARCH
process satisfying (\ref{3.0}), but with initial
values $\widetilde{X}_1=\cdots=\widetilde{X}_p=\sqrt{\omega}$.
(The latter condition is imposed to ensure that $E\widetilde{X}_1^4,\ldots,E\widetilde{X}_p^4$ are
guaranteed to be finite.)
Since $\widetilde{X}_n\stackrel{d}{\longrightarrow} X_0$ it follows from Theorem~III.6.31 in
\citet[page~58]{Pol84} that we can construct a coupling of these random variables where we have
almost sure convergence rather than convergence in probability. Hence, we obtain by Fatou's lemma that
\begin{equation}
\label{pl41.1}
EX_0^k \,\leq\, \lim\,\inf_{n\to\infty} E\widetilde{X}_n^k, \quad \mbox{ for } k=2,4.
\end{equation}
\begin{itemize}
\item[(i)]
It follows from (\ref{3.0}) that, for $t>p$,
\begin{displaymath}
E \widetilde{X}_t^2 = \widetilde{v}_t
\,\leq\, \omega \,+\, \alpha \max\{E \widetilde{X}_{t-1}^2, \ldots, E \widetilde{X}_{t-p}^2\}.
\end{displaymath}
Let $Z_t=\max\{ E \widetilde{X}_t^2, \ldots, E \widetilde{X}_{t-p+1}^2\}$.
We obtain from the previous display the recursion
\begin{displaymath}
Z_t \,\leq\, \max\{ \omega \,+\, \alpha\; Z_{t-1}, Z_{t-1} \}.
\end{displaymath}
Therefore,
\begin{displaymath}
E \widetilde{X}_t^2 \,\leq\, \frac{\omega}{1-\alpha},
\end{displaymath}
which yields in conjunction with (\ref{pl41.1}) that (i) holds true.
\item[(ii)]
If $X\sim\Skellam(v/2,v/2)$, then $EX^4=v+3v^2$.
Hence, for $t>p$
\begin{eqnarray*}
E \widetilde{X}_t^4
& = & \omega \,+\, \sum_{i=1}^p \alpha_i E \widetilde{X}_{t-i}^2
\,+\, 3\; E\left[ (\omega \,+\, \sum_{i=1}^p \alpha_i \widetilde{X}_{t-i}^2)^2 \right] \\
& \leq & \omega \,+\, 3\;\omega^2 \,+\, (1 \,+\, 6\omega) \; \alpha \; \max\{E\widetilde{X}_{t-1}^2,\ldots,E\widetilde{X}_{t-p}^2\} \\
& & {} \,+\, 3\; \alpha^2\; \max\{E\widetilde{X}_{t-1}^4,\ldots,E\widetilde{X}_{t-p}^4\}.
\end{eqnarray*}
With $\bar{Z}_t=\max\{ E \widetilde{X}_t^4, \ldots, E \widetilde{X}_{t-p+1}^4\}$ and
$\displaystyle \bar{\omega}=\omega+3\omega^2+(1+6\omega)\alpha\cdot\frac{\omega}{1-\alpha}$, we obtain the recursion
\begin{displaymath}
\bar{Z}_t \,\leq\, \max\{ \bar{\omega} \,+\, 3\; \alpha^2 \; \bar{Z}_{t-1}, \bar{Z}_{t-1} \},
\end{displaymath}
which leads to
\begin{displaymath}
E\widetilde{X}_t^4 \,\leq\, \frac{\bar{\omega}}{1 \,-\, 3\; \alpha^2}.
\end{displaymath}
\medskip

(ii) follows now from (\ref{pl41.1}).
\end{itemize}
\end{proof}
\medskip

{\lem
\label{L4.2}
Let $((\widetilde{X}_t,\widetilde{\lambda}_t))_{t\in\Z}$ be a stationary
Skellam-ARCH process satisfying (\ref{21.1a}) and (\ref{3.0}) and with $\omega>0$.
Then the matrix $\Sigma$ defined in (\ref{3.2}) is regular.
}
\medskip

\begin{proof}[Proof of Lemma~\ref{L4.2}]
We have that
\begin{displaymath}
\Sigma \,=\, E[ Z Z^T ],
\end{displaymath}
where $Z=(1,X_p^2,\ldots,X_1^2)^T$.

Assume that $\Sigma$ is singular: then there exists some
$\gamma=(\gamma_0,\ldots,\gamma_p)^T\neq 0_{p+1}$
such that
\begin{displaymath}
0 \,=\, \gamma^T \Sigma \gamma \,=\, E[(Z^T \gamma)^2],
\end{displaymath}
which implies that
\begin{displaymath}
P\left( Z^T \gamma \,=\, 0 \right) \,=\, 1.
\end{displaymath}
This means that
\begin{displaymath}
\gamma_0 \,+\, \sum_{i=1}^p X_{t-i+1}^2 \gamma_i \,=\, 0
\end{displaymath}
holds with probability~1.
Since $\gamma_1=\cdots =\gamma_p=0$ would then imply that $\gamma=0_{p+1}$,
there exists some $i_0\geq 1$ such that $\gamma_1=\cdots=\gamma_{i_0-1}=0$ and $\gamma_{i_0}\neq 0$.
Then
\begin{displaymath}
X_{t-i_0}^2 \,=\, \frac1{ \gamma_{i_0} }\left\{ \gamma_0 \,+\, \sum_{i=i_0+1}^p \gamma_i X_{t-i+1}^2 \right\}
\end{displaymath}
that is, $X_{t-i_0}^2$ is fully determined by the past values of the count process.
This, however, leads to a contradiction since
\begin{displaymath}
X_{t-i_0}\mid \F_{t-i_0-1} \,\sim\, \Skellam(v_{t-i_0}/2,v_{t-i_0}/2)
\end{displaymath}
with $v_{t-i_0}\geq\omega >0$. Hence, $\Sigma$ is a regular matrix.
\end{proof}
\bigskip

\noindent
{\bf Acknowledgement}
We are especially thankful to Bozidar Popovic and Miroslav Ristic for initiating us to investigate Skellam models.
Various preliminary discussions and suggestions concerning simple Skellam models were very fruitful and led to other considerations
beyond those in the current paper.\medskip

\noindent
This work was funded by CY Initiative of Excellence (grant "Investissements d'Avenir" ANR-16-IDEX-0008),  Project "EcoDep", PSI-AAP2020-0000000013.

%%%%%%%%%%%%%%%%%%%%%%%%%%%%%%%%%%%%%%%%%%%%%%%%%%%%%%%%%%%%%%%%%%%%%%%%%%%%%%%%%%%%%%

\bibliographystyle{harvard}

\end{document}